\documentclass[12pt, epic, eepic, reqno]{amsart}
\usepackage{amssymb}

\usepackage{epsfig}
\usepackage{psfrag}

\newtheorem{theorem}{Theorem}[section]
\newtheorem{lemma}[theorem]{Lemma}
\newtheorem{cor}[theorem]{Corollary}
\newtheorem{proposition}[theorem]{Proposition}
\newtheorem{claim}[theorem]{Claim}

\newcommand{\Z}{\mathbb Z}
\newcommand{\C}{\mathbb C}

\newcommand{\al}{\alpha}

\newcommand{\tht}{\theta}

\newcommand{\R}{\mathbb R}
\newcommand{\Q}{\mathbb Q}
\newcommand{\N}{\mathbb N}
\newcommand{\comment}[1]{}

\newcommand{\frp}{\mathfrak P}
 
\title[On Haagerup's list of potential principal graphs]{On Haagerup's list of potential principal graphs of subfactors}

\author{Marta Asaeda}

\address{Department of Mathematics, University of California Riverside,  900 Big Springs Drive, Riverside, CA, 92521,  USA} 

\email{\tt marta@math.ucr.edu}

\author{Seidai Yasuda}

\address{Research Institute for Mathematical Sciences, Kyoto University, Kitashirakawa, 
Sakyo-ku, Kyoto 606-8502, Japan} 

\email{\tt yasuda@kurims.kyoto-u.ac.jp}

\thanks{The first author was sponsored in part by NSF grant
  \#DMS-0504199.}
  
\begin{document}
 
 \begin{abstract}
%
We show that any graph, in the sequence given by Haagerup in 1991 as
that of candidates of principal graphs of subfactors, is not realized
as a principal graph except for the smallest two.  This  settles the
remaining case of a previous work of the first author.
  \end{abstract}
\maketitle
 
 \section{Introduction} 
\begin{figure}
\psfrag{Gk}{$\Gamma_k:=$}
\psfig{figure=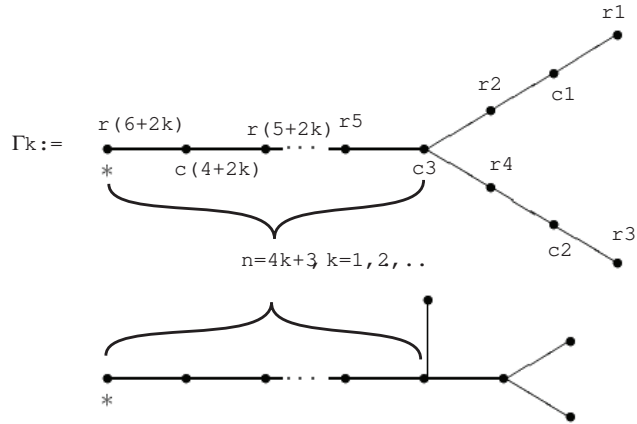,height=6cm}
\caption{The pairs of graphs (2) in the list of Haagerup}
\label{graph}
\end{figure}

This paper completes the proof that the pairs of graphs as in Fig. \ref{graph} are not realized as (dual) principal graphs of any subfactor  for $n>7$. 
 These graphs are a part of the list of graphs given by Haagerup in 1991 in \cite[\S7]{H5} as candidates which might be realized as (dual) principal graphs of subfactors. Bisch  proved that a subfactor with (dual) principal graph (4) in \cite[\S7]{H5} does {\it not} exist \cite{Bs8} by checking the inconsistency of fusion rules on the graph. Haagerup and the first author proved that two pairs of graphs:  the case 
 $n=3$ of (2) (see Figure \ref{graph})
  as well as the case (3) in \cite[\S7]{H5}, are realized as (dual) principal graphs of subfactors, and that such subfactors are unique respectively (\cite{AH}). 
 The remaining problem was whether the graphs for the case $n >3$ of (2) as in  Figure \ref{graph} would be realized as (dual) principal graphs of subfactors. Haagerup proved that the obstruction, as found for the case (4) by Bisch, does not exist on any of the pairs of the graphs in (2). Moreover, he proved that a unique biunitary connection exists for each pair of the graphs (\cite{Hp}). 
 For the case $n=7$, it was numerically checked by Ikeda that the biunitary connection should be flat (\cite{Ik}). 
 
  In 2005, Etingof, Nikshych, and Ostrik showed in \cite[Theorem 8.51]{ENO}, that the index of a subfactor has to be a cyclotomic integer, namely an algebraic integer that lies in a cyclotomic field. The result is essentially based on the result by A. Coste and T.Gannon in \cite{CG}, that shows that the entries of the $S$-matrix of a modular tensor category are in some cyclotomic field. This implies that if the square of the Perron-Frobenius eigenvalue (PFEV) of a graph is not a cyclotomic integer, the graph cannot be the (dual) principal graph of a subfactor. Utilizing this new fact, 
 the first author proved that the graphs in  Figure \ref{graph} are not (dual) principal graphs for $n=4k+3$ for $1<k \leq 27$ by showing that for each $1<k \leq 27$, the Galois group of the minimal polynomial $m_k$ of the square $d_k$ of  PFEV of each graph is not abelian: it is actually a symmetric group.  By the Kronecker-Weber theorem (\cite{Wash}), this implies that the $d_k$'s for $k$ in said range, are not cyclotomic integers. The first author also checked that for the case $k=1$, $d_1$ is a cyclotomic integer. Kondo's result in \cite{Kon}, that implies  that the Galois group of an irreducible polynomial with square-free discriminant should be symmetric, played an essential role there. 
 
In this paper we prove, by further utilizing algebraic number theory,  that none of the graphs in Figure \ref{graph} can be realized  a (dual) principal graph for $k>1$. We prove that the $d_k$'s for $k>1$ are not only not cyclotomic integers, but actually the field extension $\Q(d_k)$ over $\Q$ is not even a Galois extension: notice that if $\Q(d_k)$ was contained in some cyclotomic field, the extension $\Q(d_k)/\Q$ is necessarily Galois, since it corresponds to a subgroup of an abelian group, which is automatically a normal subgroup. 

 The first author would like to thank T.~Banica for valuable discussions, especially for bringing  \cite{BB} to attention, which contained a change of variable used in \S\ref{irred}, and D.~Bisch, V.~Jones and Y.~Kawahigashi for pointing out the result in \cite{ENO}.  M.A.~ also thanks RIMS for hospitality during the visit in May 2007, that made this collaboration possible.

\section{Essential tools from algebraic number theory} 
In the following, we list some theorems in algebraic number theory necessary for later discussion.  Most of them are directly cited from references.   We give all the proofs for the statements for which we could not find a reference. 

\begin{proposition} 
\label{unit}
 Let $\xi$ is an algebraic integer such that all the conjugates have the complex absolute value equal to one. Then $\xi$ is a root of unity. 
\end{proposition}
\noindent
{\bf Proof.} \\
Let $n$ be the number of the conjugates of $\xi$. For any $\epsilon$, there is $N$ such that 
$$ | \xi^N -1| < \epsilon/2^{n-1}. $$ 
Let $P:= \prod_{\xi'} (\xi'^N-1)$, where the product is taken over all conjugates $\xi'$'s of $\xi$. Then 
$$|P| = \prod_{\xi'} | \xi'^N-1| \leq (\prod_{\xi' \neq \xi} 2)  |\xi^N-1| < \epsilon.$$ 
Therefore we may choose $N$ so that $P$ is arbitrarily close to $0$. On the other hand, $\xi^N-1$ is also an algebraic integer, and its conjugates are given by $(\xi'^N-1)$'s. Therefore they are roots of an irreducible  monic polynomial in $\Z[x]$, thus $P \in \Z$. This means $P=0$, i.e.~ $\xi'^N-1=0$ for some $\xi'$. Then all the conjugates of $\xi'^N-1$ are also $0$, this implies $\xi^N-1=0$. Thus $\xi$ is a root of unity. \qed \\ \ \\
%
 %
 %
 %
 %
 %
 \indent
The rest of this section is devoted to a brief explanation of Hilbert's theory on ramification of ideals, which plays a key role in our argument, and to listing the theorems we use. 
 
 Let $K$ be a finite extension of $\Q$, namely a field generated by finitely many algebraic numbers. We denote by $O_K$ {\it the ring of integers of K}, namely the set of algebraic integers contained in $K$. For example, $O_\Q=\Z.$ 
 
 Let $p$ be a prime number. It generates a prime ideal $(p)$ in $\Z$. Now, consider the ideal  $pO_K$, generated by $p$ in $O_K$. This is not generally a prime ideal. Since $O_K$ is a {\it Dedekind domain} (\cite{lecturenote}, 3.1), it factorizes into a product of prime ideals uniquely: 
 $$ pO_K =\frp_1 ^{e_1} \cdots \frp_g^{e_g},$$
 where $\frp_i$'s are distinct prime ideals of $O_K$. It is easy to see that $\frp_i \cap \Z=(p)$ for all $i$. We call $e_i$ {\it the ramification index of $\frp_i$}. For a prime ideal $\frp$ of $O_K$, $O_K/\frp$ is a field. Consider the composition of the maps
 $$ \Z \stackrel{\iota}{\hookrightarrow} O_K \stackrel{\pi}{\twoheadrightarrow} O_K/\frp.$$
 Then Ker$\pi \circ \iota=\Z \cap \frp =(p)$. Thus $\pi \circ \iota$ induces a field extension $k:=\Z/p\Z \hookrightarrow O_K/\frp$. We call $[O_K/\frp, k]=:h(\frp)$ {\it the degree of $\frp$ over $k$.} 
 
 The ramification theory concerns the factorization described above, for a given prime $p$ and a field extension $K$. There is the following beautiful theorem. 
 
 \begin{theorem}(Dedekind, \cite{Nar}, Theorem 4.33)
 Let $d$ be an algebraic integer, $K=\Q(d)$, and $f(x) \in \Z[x]$ be the minimal polynomial of $d$. Let $p$ be a prime number that does not factor $D_f /D_K$, where $D_f$ is the discriminant of $f$, $D_K$ is the discriminant of $K$, and let $k=\Z/p\Z$. Suppose the factorization of $f$ mod $p$ is given by 
 $${\bar f}(x) \equiv {\bar f}_1^{e_1} \cdots {\bar f}_g^{e_g} \ \ {\rm mod} \ p, $$
 where ${\bar f}_i$'s are irreducible polynomials in $k[x]$. Then we have
 $$pO_K= \frp_1 ^{e_1} \cdots \frp_g^{e_g},$$
 where $\frp_i$'s are distinct prime ideals of $O_K$, and $h(\frp_i)=\deg {\bar f}_i$. 
 \end{theorem}
 Here we do not give definitions for the discriminant of a polynomial nor the discriminant of a field. In fact we do not want to deal with the discriminants, thus we need to modify this theorem for our use. We will also combine it  with the following nice theorem: 
 \begin{theorem} (\cite{Nar}, Theorem 4.6)
 Suppose $K/\Q$ is a Galois extension of degree $n$. Then for a prime $p$ we have
 $$pO_K=(\frp_1 \cdots \frp_g)^e,$$
 where $\frp_i$'s are distinct prime ideals of $O_K$, and $h(\frp_i)=h$ for all $i$ for some $h$, and we have $n=ehg.$
 \end{theorem}
 We obtain the following theorem for our use. 
 \begin{theorem}
 \label{factorpoly}
 Let $d$ be an algebraic integer, $K=\Q(d)$, and $f(x) \in \Z[x]$ be the minimal polynomial of $d$ with degree $n$. Suppose that $K/\Q$ is Galois. Let $p$ be  a prime number, and $k:=\Z/pZ$. Let  $e$, $f$, and $g$ be integers such that
 $$pO_K=(\frp_1 \cdots \frp_g)^e,$$
 where $\frp_i$'s are distinct prime ideals of $O_K$, and $h(\frp_i)=h$ for all $i=1, \dots, g.$ Then $f(x)$ factorizes mod $p$ as follows: 
 $$ {\bar f}(x)= ({  f}_1 \cdots {  f}_g)^e \ \  {\rm mod} \  p, $$
 where $f_i \in k[x]$ with $\deg {f}_i=h$ for all $i$ and each ${f}_i$ is of the form ${ f}_i=g_i^{e'_i}$, where $g_i \in k[x]$ is irreducible. 
 \end{theorem}
 \noindent
{\bf Proof.} \\ 
Let $G:={\rm Gal}(K/\Q)$  and $k_i:=O_K/\frp_i.$ Note that for $\sigma \in G$, $\sigma(\frp_i)$ is a prime ideal, and it coincides with some $\frp_j$, since $\sigma(\frp_i) \cap \Z=(p).$ 
For each $\frp_i$ we define
$$H_i:=\{ \sigma \in G | \sigma(\frp_i)=\frp_i\}.$$
Then $H_i$ is a subgroup of $G$. Consider the following surjection 
$$\psi_i: H_i \twoheadrightarrow {\rm Gal}(k_i/k)=:G_i.$$ 
Let $I_i:={\rm Ker} \ \psi_i =\{ \pi  \in H_i | \pi(a) \equiv a \ {\rm mod} \ \frp_i, \forall a \in O_K\}$. This is a normal subgroup of $H_i$, and we have $H_i/I_i \cong G_i.$ \footnote{$H_i$ and $I_i$ are called {\it decomposition group} and {\it inertia group} of $\frp_i$ respectively, see \cite{Nar}, p263.}
For each $i$, let $\sigma_i \in G$ to be so that $\sigma_i(\frp_1)=\frp_i$. Then we obtain a coset decomposition  $G= \sigma_1 H_1 \sqcup \cdots \sqcup \sigma_g H_1.$ Observe that $H_i=\sigma_i H_1 \sigma_i^{-1}, $ thus 
$G= H_1 \sigma_1 \sqcup \cdots H_g \sigma_g.$ 
Note also that $|H_i|=|H_1|$ and $ |G_i|=[k_i:k]=h$ for all $i$, thus we have $|I_i|=|I_1|$ as well. Noting that $n=|H_1|g=|I_1|hg$, we get $|I_i|=e.$ 

In $O_K[x]$ we have
$$f(x)=\prod_{\sigma \in G} (x-\sigma(d))=\prod_i \prod_{\sigma \in H_i \sigma_i} (x-\sigma(d)).$$
Let $v_i(x):= \prod_{\sigma \in H_i \sigma_i} (x-\sigma(d)) \in O_K[x]$. Since $H_i$ preserves $\frp_i$, we have $\sigma(\sigma_i(d)) \equiv \psi(\sigma) (d_i) \ \ {\rm mod} \ \frp_i$, where $\sigma \in H_i$ and  $d_i$ is the image of $\sigma_i(d)$ in $O_K/\frp_i.$   
Noting  $H_i/I_i \cong G_i$, we have
\begin{eqnarray*} 
v_i(x) &\equiv& \prod_{\sigma   \in H_i  } (x- \psi(\sigma)(d_i)) \ \ {\rm mod}  \ \frp_i \\
&\equiv& \prod_{\tau \in I_i} \prod_{\rho \in G_i} (x- \rho \psi(\tau)(d_i)) \ \ {\rm mod}  \ \frp_i \\
&\equiv& \prod_{\tau \in I_i} \prod_{\rho \in G_i} (x- \rho(d_i)) \ \ {\rm mod}  \ \frp_i  \\
&=& (\prod_{\rho \in G_i} (x- \rho(d_i)))^e \ \ {\rm mod}  \ \frp_i  
\end{eqnarray*}
Note that $f_i(x):=\prod_{\rho \in G_i} (x- \rho(d_i)) \in k[x]$, and $\deg f_i=|G_i|=h.$ Thus $v_i(x) \ \ {\rm mod} \ \frp_i \in k[x]$ as well. 

The polynomial $f_i(x)$'s may or may not be irreducible in $k[x]$. Let $F_i:=\{ \tau \in G_i | \tau(d_i)=d_i\}.$ Since $G_i$ is abelian, $F_i$ is a normal subgroup of $G_i$. Thus we have
\begin{eqnarray*}
f_i(x)&=& \prod_{\rho \in G_i/F_i} \prod_{\tau \in F_i}  (x- \rho\tau(d_i)) \\
&=& (\prod_{\rho \in G_i/F_i} (x-\rho(d_i)))^{e'_i}, 
\end{eqnarray*}
where $e'_i=|F_i|.$ Since $\rho(d_i)$ runs through all the conjugates of $d_i$, $g_i(x):=\prod_{\rho \in G_i/F_i} (x-\rho(d_i))$ is the minimal polynomial of $d_i$ and $g_i(x) \in k[x]$. 

Now, since $g_i(x)$ mod $\frp_i$ $\in k[x]$, we have $v_i(x)=g(x)^{e'_i e}$ mod $p$. Altogether we have desired factorization of $f(x)$ mod $p$, 
$${\bar f}(x)= v_1(x) \cdots v_g(x) = (f_1 \cdots f_g)^e \ \ {\rm mod} \ p,$$
where for each $i$ $\deg f_i=h$, $f_i =g_i^{e'_i}$, and $g_i$ is irreducible. 
\qed 
\section{Minimal polynomials}
\label{irred}
Let $d_k$ be the square of PFEV of the graph $\Gamma_k$ in Fig. \ref{graph}. In \cite{A2}   the adjacency matrix $A_k$ of $\Gamma_k$ was given, which is of the size  $(4+2k) \times (6+2k)$. The characteristic polynomial of the matrix $N_k:={A_k}^t A_k$ divided by $(x-2)^2$, which is denoted by $q_k(x)$, satisfies the following recursive formula
\begin{eqnarray*}
\label{recursive}
q_k(x)&=&(x^2-4x+2)q_{k-1}(x)-q_{k-2}, \\
q_0(x)&=&x^2-5x+3, \\
 q_1(x)&=&(x^3-8x^2+17x-5)(x-1).
 \end{eqnarray*}
and thus computed as follows:
$$q_k(x)=A(x) a(x)^{2k} + B(x) b(x)^{2k},$$
where $a(x) = (2 - x + \sqrt{x^2 - 4x})/2$, $b(x) = (2 - x -\sqrt{x^2 - 4x})/2$, $A(x)=\frac{-1}{a(x)^2-b(x)^2}(q_0(x)b(x)^2-q_1(x))$, and $B(x)=\frac{1}{a(x)^2-b(x)^2}(q_0(x)a(x)^2-q_1(x))$. The largest root of $q_k$ is $d_k$.

 In this section we prove the following theorem conjectured in \cite{A2}. 
 \begin{theorem}
 \label{irred}
Let 
$$
r_k(x)=  \begin{cases} 
q_k(x)/(x-1), &   {\rm if} \; k \equiv 1 \ {\rm mod}\  3, \\
q_k(x), &  {\rm else}.  
\end{cases}
$$
Then $r_k(x)$ is irreducible for any $k$, thus it is the minimal polynomial of $d_k$. 
\end{theorem}

\begin{figure}[h]
\psfig{figure=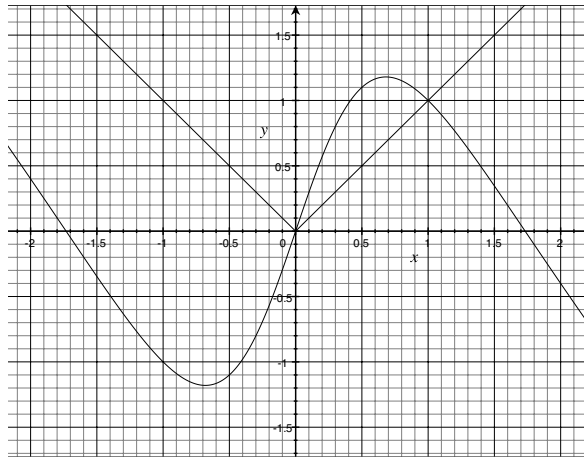,height=6cm}
\label{f(x)|x|}
\caption{The graphs for $y=|x|$ and $y=f(x)=\frac{3x-x^3}{1+x^2}$}
\end{figure}

One immediately sees that the polynomials $q_k(x)$'s are ugly: indeed
\begin{eqnarray*}
q_2(x) &=& x^6-13x^5+63 x^4-140 x^3 +142 x^2 -59 x +7, \\
q_3(x) &=& x^8-17x^7 +117 x^6-418 x^5+827 x^4 -898 x^3 +502 x^2 -124 x+9, 
\end{eqnarray*}
and so on. It is hard to see any pattern as $k$ varies.  However, by the change of variable used in \cite{BB}, we obtain better polynomials.  We define
$$P_k(q):= q_k(x)|_{x=q+q^{-1}+2} q^{2k+2}. $$
 The polynomials $P_k$'s satisfy the recursive formula
 \begin{eqnarray}
\label{recursive}
P_k(q)&=&(q^4+1)P_{k-1}(x)-q^4 P_{k-2}, \\
P_0(q)&=&q^4-q^3-q^2-q+1, \\
 P_1(q)&=&q^8-q^7-q^6-q^5+q^4-q^3-q^2-q+1.
 \end{eqnarray}
 Thus we obtain 
 $$P_{k-1}(q)=q^{4k}-q^{4k-1}-q^{4k-2}-q^{4k-3}+q^{4k-4}- \cdots -q^5+q^4-q^3-q^2-q+1. $$
 for any $k \geq 1$. 
 Our goal is to prove the following theorem, which is stronger than Theorem \ref{irred}.
 \begin{theorem}
 \label{irredq}
 For each $k  \geq 1$, let 
 $$ R_{k-1}(q):=
 \begin{cases}
P_{k-1}(q) \ {\text  if } \  k \neq 2 \ \mbox{{\rm mod} } 3 \\
P_{k-1}(q)/ (q^2+q+1) \  {\text{if} } \   k =2 \ \mbox{{\rm mod} } 3. 
 \end{cases}
 $$
Then $R_{k-1}(q)$ is irreducible. 
 \end{theorem}
 \begin{proposition}
 \label{roots}
 Let $k \geq 0$.  
 \begin{itemize}
 \item[(1)] Then there exists unique $\al \in (0,1)$ such that $P_k(\al)=0.$ \hfill\break 
 ( $\Leftrightarrow$ (1)' there exists unique $\al' >1$ such that $P_k(\al')=0$.) 
 \item[(2)] If $\beta \in \C$ is a root of $P_k$, then $\beta = \al, \al'$, or $| \beta|=1$. 
 \end{itemize}
 \end{proposition}
 This, together with Proposition \ref{unit}, implies the following: 
 \begin{cor}
 \label{factorize}
Suppose $P_k$ factorizes into the product of irreducible polynomials as follows: 
$$P_k(q)=P_{k,1}(q) \dots P_{k, r}(q),$$
 and suppose $P_{k,1}(\al)=0$. Then $P_{k,1}(\al')=0$, and for $i \geq 2$, all the roots of $P_{k,i}$ are roots of unity. 
 \end{cor}
 \noindent
 {\bf Proof of Proposition \ref{roots}.} \\
%
(1): Notice that $P_k(0)=1>0$, $P_k(1)=-2k-1<0$, thus there exist a root $\alpha$ of $P_k$ in $(0,1)$. We show that it is unique. It suffices to show that $P'_k <0$ on $(0,1)$. For $k=0$, $P_0(q)=q^4-q^3-q^2-q+1$, so $P'_0(q)=4q^3-3q^2-2q-1$. 
Since $q^3 < q^2 < q$ on $(0, 1)$,
$$P'_0(q) < (3q^3+q)-3q^2-2q-1=-q-1 < 0$$
holds  in $(0,1)$.  For general $k$, since 
$$P_{k-1}(q)-P_{k-2}(q) =(q^{4k}-q^{4k-1}-q^{4k-2}-q^{4k-3})=q^{4k-3}(q^3-q^2-q-1),$$
we have
$$(P_{k-1}(q)-P_{k-2}(q))'=(4k-3)q^{4k-4} (q^3-q^2-q-1) + q^{4k-3} (3q^2-2q -1).$$
It is easily checked that $(q^3-q^2-q-1),  (3q^2-2q -1) <0$ in $(0,1)$. Thus $P'_k(q) < P'_{k-1}(q)< ... < P'_0(q)<0$ in $(0,1)$. (1)' is immediate from the fact that $P_k(q^{-1})q^{4(k+1)} = P_k(q)$. 
 \\ \ \\
\noindent
 (2): 
Notice that $q^4-q^3-q^2-q \geq 0$ for $q  \leq 0.$  Therefore 
$$P_{k-1}(q)=\sum_{l=1}^k (q^4-q^3-q^2-q) q^{4l-4} +1 > 0 $$
for $q  \leq 0$, which implies that $P_{k-1}(q)$ has no non-positive real root. Thus the only real roots of $P_{k-1}(q)$ are $\alpha$ and $1/\alpha$. On the other hand, recall that the matrix  $N_k:={A_k}^t A_k$ is symmetric, thus all the eigenvalues are real. Therefore all the roots of $q_{k-1}(x)$ are real. If $\beta$ is a root of $P_{k-1}(q)$, $\beta +1/\beta=r$ is a root of $q_{k-1}(x)$, which is real, and $\beta$ is a root of $t^2-rt+1=0$. This implies that $\beta$ is real or $|\beta|=1$. \qed \\ \ \\
%
%
{\bf Proof of Theorem \ref{irredq}.}\\ 
For $k \neq 2$ mod $3$, we show that $P_{k-1}(q)$ is irreducible. From Cor. \ref{factorize}, it suffices to show that $P_{k-1}(q)$ has no root which is a root of unity. 
 Let 
$$Q_{k-1}(q):=P_k(q)(q^4-1)=q^{4k+4}-q^{4k+3}-q^{4k+2}-q^{4k+1}+q^3+q^2+q-1.$$
Note that the roots of $Q_{k-1}(q)$ are the roots of $P_{k-1}(q)$ except for $q= \pm 1, \pm i$: it is easy to check that they are not roots of $P_{k-1}(q)$.  Thus it suffices to show that $Q_{k-1}(q)$ has no root which is a root of unity except for those. 
Let $\beta=e^{2 \pi i \tht}$, where $\theta \in [0,1)$, and suppose $Q_{k-1}(\beta)=0$.
 Notice that 
 \begin{eqnarray*}
  Q_{k-1}(q)= && q^{2k+2}((q^{2k+2}-q^{-(2k+2)})-(q^{2k+1}-q^{-(2k+1)}) \\
 && -(q^{2k}-q^{-2k})-(q^{2k-1}-q^{-(2k-1)})).
 \end{eqnarray*}
   Thus $Q_{k-1}(\beta)=0 \Leftrightarrow$
   \begin{eqnarray*}
&& \sin 2(2k+2)\pi \tht -\sin 2(2k+1)\pi \tht -\sin 4k \pi \tht -\sin2(2k-1)\pi \tht \\
 &=& 2 \sin 2 (2k+\frac{1}{2}) \pi \tht \cos 3 \pi \tht -2 \cos 2 (2k+\frac{1}{2}) \pi \tht \sin \pi \tht =0
 \end{eqnarray*} 
 $\Leftrightarrow   \tht= \frac{1}{2}$ or
\begin{equation}\label{flat}\tag{$\flat$} 
 \tan(4k+1)\pi \tht =\frac{\sin 3\pi \tht}{\cos \pi \tht}. 
\end{equation}
Notice that $$\frac{\sin 3\pi \tht}{\cos \pi \tht} =\frac{3 \tan \pi \tht -\tan^3 \pi \tht}{1+\tan^2 \pi \tht}.$$ Therefore, 
\begin{equation}\label{sharp}\tag{$\sharp$} 
\eqref{flat} \Leftrightarrow  \tan(4k+1)\pi \tht =f(\tan \pi \tht), 
\end{equation}
where $f(x):=\frac{3x - x^3}{1+x^2}$. 
Thus we need to show that there is no $\tht \in \Q \cap [0,1)$ satisfying the equation \eqref{sharp}
except for $\tht=\frac{1}{4}, \frac{3}{4}$ and $0$. Similarly, for $k = 2 \mod 3$, we need to show that the only roots of $Q_{k-1}(q)$ which are roots of unity are  the roots of $(q^4-1)(q^2+q+1)$. So we need to show that there is  no $\tht \in \Q \cap [0,1)$ satisfying the equation  \eqref{sharp} except for $\tht= \frac{1}{3}, \frac{2}{3}$ in addition. 
\\ \ \\ \indent
Suppose there is $\tht =  \frac{m}{N} \in  [0,1)$  satisfying  \eqref{sharp}, where $m, N \in \N$, $N \geq 3$, and  $(m,N)=1$. 
\begin{lemma}
\label{Nx}
For $ \forall b \in (\Z/N\Z)^\times$, $b \tht$ satisfies  \eqref{sharp}. 
\end{lemma}
\noindent
{\bf Proof.} \\
Let $K$ be the splitting field of $Q_{k-1}(q)$ and $G={\rm Gal}(K/\Q)$. By the assumption $e^{2\pi i \tht} \in K$, thus $K \supset \Q(e^{\frac{2\pi i}{N}})$. Observe ${\rm Gal}(\Q(e^{\frac{2\pi i}{N}})/\Q)= (\Z/N\Z)^\times$, where the action of $b\in  (\Z/N\Z)^\times$ is given by \\
$\sigma_b \in {\rm Gal}(\Q(e^{\frac{2\pi i}{N}})/\Q)$, $\sigma_b (e^{\frac{2\pi i}{N}})=e^{\frac{2\pi ib}{N}}$. Take $g \in G$ such that ${\bar g} =\sigma_b \in G/{\rm Gal}(K/\Q(e^{\frac{2\pi i}{N}}))$, then $g (e^{\frac{2\pi im }{N}})= \sigma_b(e^{\frac{2\pi im }{N}})= e^{\frac{2\pi imb }{N}}$, thus $e^{\frac{2\pi imb }{N}}$ is a root of $Q_{k-1}(q)$ as well, thus $\frac{mb}{N}=b \tht$ satisfies  \eqref{sharp}.  \qed 
\\ \ \\ \indent
Therefore, without loss of generality we choose $\tht$ so that $|\frac{1}{2} - \tht|$ will be the minimum among the choices of $\tht$, which implies that $|\tan \pi \tht|$ is the maximum. We may choose so that $\frac{1}{2} - \tht>0$, thus $\tan \pi \tht>0$. More specifically, we choose 
$$
\tht=
\begin{cases}
\frac{N-1}{2N} \mbox{ {\text{if}  {\it N} is odd} } \\
\frac{\frac{N}{2}-2}{N}  \mbox{ {\text{if} } } N\equiv 2 \mod 4 \\
\frac{\frac{N}{2}-1}{N}  \mbox{ {\text{if} } } N\equiv 0 \mod 4 
\end{cases}
$$
\begin{lemma}
$$ {\rm gcd} (N, 4k+1)= 1 \mbox{ {\text or} }3. $$
In particular, for $k \neq 2 \mod 3$, $ {\rm gcd} (N, 4k+1)= 1$. 
\end{lemma}
\noindent
{\bf Proof.}  \\
Let $d:=\gcd(N, 4k+1)$, and 
$$S:=\{ b \in (\Z/N\Z)^\times | \tan (4k+1) \pi b \tht = \tan (4k+1) \pi \tht \}. $$ 
Then $b \in S \Leftrightarrow (4k+1)b=(4k+1) \mod N \Leftrightarrow b =1 \mod \frac{N}{d}.$
\begin{lemma}
\label{upperbound-d}
$$|S| \geq \varphi(d)=:|(\Z/d\Z)^\times|.$$
\end{lemma}
We prove this lemma later on. Using this lemma will give an upperbound of $d$. 
For $b \in S$, we have 
$$ f(\tan \pi b \tht)=\tan (4k+1) b \pi \tht =\tan (4k+1) \pi \tht.$$ 
Note that the last term is fixed. Since $\deg f=3$, there are at most three solutions to $f(x)=const$. Therefore we obtain $3 \geq |S| \geq \varphi(d)$. Noting that $d | 4k+1$, $d$ needs to be odd. Thus we get $d=1$ or $3$.  $d=3$ is possible only if $3 | 4k+1 \Leftrightarrow k=2 \mod 3$. \qed
\\ \ \\ 
\noindent
{\bf Proof of Lemma \ref{upperbound-d}.} \\
There is a natural group homomorphism 
$$\psi: (\Z/N\Z)^\times  \longrightarrow (\Z/(N/d)\Z)^\times.$$
Observe that $\ker \psi=S$. Thus 
\begin{equation}\label{star}\tag{$\star$}
\varphi(N)/|S| \leq \varphi(N/d).  
\end{equation}
There is a formula for computing $\varphi$ (\cite{eulerwiki}): for $n=p_1^{e_1} \cdots p_r^{e_r}$, where $p_i$'s are distinct primes, we have
$$\varphi(n)=(1-\frac{1}{p_1}) \cdots (1-\frac{1}{p_r}) \cdot n.$$
Applying this formula to $\eqref{star}$ we obtain $|S| \geq \varphi(d)$. \qed
\\ \ \\
\indent
We return to the proof for Theorem \ref{irredq}. \\ \ \\
\noindent
{\bf Case 1}: $k \neq 2 \mod 3$. In this case $d:= {\rm gcd} (N, 4k+1)=1$. Let $\tht'$ to be so that $(4k+1)\tht'=\tht$. Since $4k+1 \in (\Z/N\Z)^\times$, $\tht'$ satisfies  \eqref{sharp}. Then 
\begin{equation}\label{natural}\tag{$\natural$}
|\tan \pi \tht'| \leq |\tan \pi \tht| = |\tan \pi (4k+1)\tht'| =|f(\tan \pi \tht')|. 
\end{equation}
We find the range of $x$ so that $|x| \leq |f(x)|.$ The graphs of $y=|x|$ and $y=f(x)$ is given in Fig. \ref{f(x)|x|}. 
For $x \geq 0$, since $f(x)=x \Leftrightarrow \frac{4x}{1+x^2}=2x \Leftrightarrow x=0$ or $\pm 1 $.  It is easy to check that $f(x) \geq x$ for $0 \leq x \leq 1$, and $f(x) <x$ for $x>1$. Since $f(x)$ is an odd function, we have $|x| \leq |f(x)|$ for $|x| \leq 1$. Since $f(x)=-x \Leftrightarrow \frac{4x}{1+x^2} =0 \Leftrightarrow x=0$, we have no other range of $x$ satisfying $|x| \leq |f(x)|$. This, together with \eqref{natural}, implies   $|\tan \pi \tht'| \leq 1$.

We shall find the maximum of $f(x)$ in  $|x| \leq 1$. 
\begin{eqnarray*}
f'(x) &=& \frac{(3-3x^2)(1+x^2)-(3x-x^3) 2x}{(1+x^2)^2} =\frac{3-6x^2-x^4}{(1+x^2)^2} =0\\
& \Leftrightarrow& 3-6x^2-x^4=0 \Leftrightarrow x^2=-3 \pm 2 \sqrt{3} 
\end{eqnarray*}
Thus the critical points are given by $x=\pm \sqrt{2\sqrt{3}-3} <1$. One may easily check that this gives local maxima for $|f(x)|$, with the value $f(\sqrt{2\sqrt{3}-3}) =: \gamma \approx 1.17996.$ Thus $|f(x)| \leq  \gamma$ for $|x| \leq 1$. 

From \eqref{natural}, we have $   |\tan \pi \tht| =  |f(\tan \pi \tht')| \leq \gamma.$ Recall that $\tht$ was given explicitly for each $N$. We now examine each case. 
\begin{itemize}
\item $N=3$: $\tht =\frac{1}{3}. \tan \pi/3=\sqrt{3} > \gamma$. Since $ \frac{1}{2} > \frac{N-1}{2N} > \frac{1}{3}$ for all $N>3$, we have $\tan \pi \tht > \sqrt{3} > \gamma$ for all the odd integer  $N>3$.  
\item $N=4$: $\tht=\frac{1}{4}. \tan \pi/4 =1 < \gamma$.  
\item $N=6$: $\tht =\frac{1}{6}. \tan \pi/6=0.57\dots <\gamma$. 
\item $N=8$: $\tht=\frac{3}{8}> \frac{1}{3}.$ We have $\tan \pi \tht > \sqrt{3} >\gamma$ for all $N>8$, $N =0 \mod 4$. 
\item $N=10$: $\tht=\frac{3}{10}.  \tan \pi \frac{3}{10} = 1.37\dots > \gamma$. We have $\tan \pi \tht   > 1.18$ for all $N>10$, $N =2\mod4$. 
\end{itemize}
We need to check that $\tht=\frac{1}{6}$ is not a solution for  \eqref{sharp}. Since $4k+1 \in (\Z/6\Z)^\times$, $(4k+1) \frac{1}{6}=\frac{1}{6}\  {\text or} \ \frac{5}{6}.$ Thus 
$$\tan  (4k+1) \pi \tht = \pm \frac{1}{\sqrt{3}}.$$ 
On the other hand,  
$$f(\tan \pi \tht)=\frac{\sin \frac{ 3\pi}{6}}{\cos\frac{\pi}{6}}=\frac{1}{\sqrt{3}/2}=\frac{2}{\sqrt{3}} \neq \pm \frac{1}{\sqrt{3}}. $$
Therefore, the only rational solutions for   \eqref{sharp} in $[0,1)$ are $\tht=\frac{1}{4}, \frac{3}{4}$. Thus the polynomial $P_{k-1}(q)$ is irreducible in this case. 
\\ \ \\
{\bf Case 2}: $k=2\mod 3$. In this case $d$ can be either $1$ or $3$. Note that $3 | 4k+1$. For the case $d=1$, $N$ cannot be divisible by $3$. The proof proceeds exactly the same as for Case 1, except that we do not have to worry about $N=6$ at the end. 

For $d=3$, we have $|S| \geq 2$ from Lemma \ref{upperbound-d}. Note that for $b\in S$, $b\tht$ is a solution for  \eqref{sharp}, by Lemma \ref{Nx}. Thus 
$$b \in S \Rightarrow \tan (4k+1) \pi \tht = \tan(4k+1) b \pi \tht =f(\tan b\pi \tht).$$
Since distinct values of $b \in S$ give distinct values for $\tan b\pi \tht$, $|S| \geq 2$ implies that 
\begin{equation}
\label{*}\tag{$*$}
\tan (4k+1) \pi \tht =f(x) 
\end{equation}
has at least two solutions, and they are in the range of $|x| \leq \kappa \approx 2.542\dots, $ where $f(\kappa)=-\gamma$. Taking $b =1$, we have $\tan \pi \tht \leq \kappa$. Noting that $3 |N$, we examine each $N=3, 6, 9...$.
\begin{itemize}
\item $N=3$: $\tht =\frac{1}{3}. \tan \pi/3=\sqrt{3}  < \kappa. $
\item $N=6$: $\tht=\frac{1}{6}. \tan \frac{\pi}{6}=0.57\dots < \kappa.$ 
\item $N=9$: $\tht=\frac{4}{9}. \tan \pi \frac{4}{9}=5.67\dots >\kappa.$ Thus for odd $N>9$, $\tan \pi \tht >\kappa$. 
\item $N=12$: $\tht=\frac{5}{12}. \tan \pi  \frac{5}{12}=3.73\dots >\kappa.$ Thus for $12<N=0\mod4$, $\tan \pi \tht >\kappa. $
\item $N=18$: $\tht=\frac{7}{18}>\frac{5}{12}.$  Thus for $18 \leq N=2\mod4$, $\tan \pi \tht >\kappa. $
\end{itemize}
We check if the surviving values $\tht=\frac{1}{3}, \frac{1}{6}$ would give solutions to  \eqref{sharp}. 
For $ \tht=\frac{1}{3}$, $\tan \pi (4k+1) \frac{1}{3}=0$ since $3|4k+1$. On the other hand, 
$$f(\tan \pi\tht)|_{\tht=1/3}=\frac{\sin \pi}{\cos \frac{\pi}{3}}=0=\tan (4k+1) \frac{1}{3}.$$
Thus $\tht=\frac{1}{3}$ and $\frac{2}{3}$ are solutions. For $\tht=\frac{1}{6}$, noting that $2 \mathrel{\not|} 4k+1$, 
$3 | 4k+1$, $\tan (4k+1) \pi \tht$ is undefined. On the other hand $f(\tan \frac{\pi}{6})= \frac{2}{\sqrt{3}} $, thus  \eqref{sharp} fails. Altogether, for $k=2\mod 3$, the only solutions for  \eqref{sharp} are $\tht=\frac{1}{4}, \frac{3}{4}, \frac{1}{3},$ and $ \frac{2}{3}. $ This complete the proof of Theorem \ref{irredq}, and thus that of Theorem \ref{irred}. \qed

\section{Factorization of minimal polynomials over primes and non-cyclotomicity of $d_k$}
In this section we show that $d_k$'s are not cyclotomic integers for $k \geq 2$, which implies that the graphs $\Gamma_k$ in Fig. \ref{graph} are not principal graphs for subfactors for $k \geq 2$, which was conjectured in \cite{A2}. 

For simplicity, we prove the equivalent statement that $e_k=d_k-2$ is not cyclotomic integers for $k \geq 2$. We shift the variable of all the polynomials accordingly: 
\begin{itemize}
\item The minimal polynomial for $e_k$ is  $m_k(x):=r_k(x+2)$. 
\item $p_k(x):=q_k(x+2)$. 
\end{itemize}
 Then 
$$
p_{k-1}(x)=
\begin{cases}
m_{k-1} (x) & \text{if} \  k\neq 2\mod 3, \\
(x+1) m_{k-1}(x) & \text{if} \  k= 2\mod 3. 
\end{cases}
$$
It relates to $P_k(q)$ by $p_{k-1}(q+q^{-1})q^{2k}=P_{k-1}(q).$  The polynomial $p_k(x)$ satisfies the recursive formula: 
\begin{eqnarray*}
\label{recursivep}
p_k(x)&=&(x^2-2)p_{k-1}(x)-p_{k-2}, \\
p_0(x)&=&x^2-x-3, \\
 p_1(x)&=&(x^3-2x^2-3x+5)(x+1).
 \end{eqnarray*}

In the rest of this section we show the following theorem: 
\begin{theorem}
\label{notgalois}
The field extension $\Q(e_{k-1})/\Q$ is not Galois for $k \geq 3.$ Thus the graphs $\Gamma_k$ in Fig.1 are not principal graphs of subfactors for $k \geq 2.$
\end{theorem}
 \noindent
 {\bf Proof.} \\
 Let $k \geq 3$ for the rest of this section.  
Suppose  $\Q(e_{k-1})/\Q$  was a Galois extension. It coincides with the splitting field of the minimal polynomial $m_{k-1}(x)$ of $e_k$. We use Theorem \ref{factorpoly} to derive a contradiction. First we look for a suitable prime number.   The following is obtained by easy computations using the recursive formula. 
\begin{claim}
\label{zero}
 \begin{eqnarray*}
 p_{k-1}(0)&=&(-1)^k (2k+1), \\
  p'_{k-1}(0)&=&(-1)^k k.
  \end{eqnarray*}
 \end{claim}
 This implies the following. 
 \begin{proposition}
 \label{factor}
 Suppose  $\Q(e_{k-1})/\Q$ is a Galois extension of $\Q$. Then for a prime $p$ such that $p | 2k+1$, $p \mathrel{\not|} k$, we have 
 $$m_{k-1}(x)= x \prod_{0 \neq a \in \Z/p\Z}  (x-a)^{n_a}\mod p. $$
 \end{proposition}
 Note that the condition $p \mathrel{\not|} k$ is obviously redundant, and that $x \mathrel{\not|} x+1\mod p$. \\  \ \\
 \noindent
 {\bf Proof.} \\
 Claim \ref{zero} implies that $x | p_{k-1} \mod p$, $x^2 {\not|}  \ p_{k-1} \mod p$. Thus, in the setting of Theorem \ref{factorpoly} we have $e=h=1$, therefore $m_{k-1}(x) \mod p$  factorizes into a product of  linear terms.  \qed \\ 
 
 In the following, we find a suitable  prime $p$ to derive a contradiction to the above proposition. 
 
   \begin{lemma}
 \label{multi}
 If $p>3$, $n_a \leq 4.$
 \end{lemma}
 \noindent
{\bf Proof.} \\
 Consider the fourth derivative of $Q_{k-1}(q)=P_k(q)(q^4-1)=q^{4k+4}-q^{4k+3}-q^{4k+2}-q^{4k+1}+q^3+q^2+q-1$:
 \begin{eqnarray*}
 Q^{(4)}_{k-1}(q)&=&(4k+4)(4k+3)(4k+2)(4k+1)q^{4k}\\
 &-& (4k+3)(4k+2)(4k+1)4k q^{4k-1} \\
 &-& (4k+2)(4k+1)4k(4k-1) q^{4k-2} \\
 &-& (4k+1)4k(4k-1)(4k-2) q^{4k-3}. 
 \end{eqnarray*}
Let $p | 2k+1$, $p>3$. Then 
 $$ Q^{(4)}_{k-1}(q) \equiv -(4k+1)4k(4k-1)(4k-2) q^{4k-3} \not\equiv 0 \mod p.$$
Thus, for $\beta$ in an algebraic closure of $\Z/p\Z$, $ Q^{(4)}_{k-1}(\beta) \equiv 0 \mod p$ only if $\beta=0$. Note that $q=0$ is not a root of $Q_{k-1}(q)$. This implies that the multiplicities of roots of $Q_{k-1}(q) \mod p$ cannot be more than four, nor can the multiplicities of the roots of $P_{k-1}(q) \mod p$. Recall that $p_{k-1}(q+q^{-1})q^{2k}=P_{k-1}(q)$. There is a one to one correspondence between factors $(x-a) \Leftrightarrow (q^2-aq+1)$.  Therefore $n_a \leq 4.$ 
 \qed \\ \ \\
 \indent
 In the following, there is a slight difference in arguments for $k\neq 2 \mod 3$ and $k=2 \mod 3$. We deal with each case one by one. 
 \subsection{The case $k\neq 2 \mod 3$} 
 \label{knequiv2}
 \ \\
 
 \noindent
 %
{\bf Case 1}: $2k+1$ is not a prime, nor a power of $3$. 
\\
 By the assumption, there is a prime number $p \neq 2k+1, 3$ that divides $2k+1$. Since $2 \mathrel{\not|} 2k+1$, $2k+1$ is divisible by some number larger or equal to $5$, thus  $p \leq \lfloor \frac{2k+1}{5} \rfloor$, where  by  $\lfloor c \rfloor$ for $c \in \R$ we denote the largest integer dominated by $c$.

 Suppose that $\Q(e_{k-1})/\Q$ is Galois. By Proposition \ref{factor} and  Lemma \ref{multi}, and that $\deg p_{k-1}=2k$, we need at least $ \lceil \frac{2k-1}{4} \rceil +1$ distinct elements in $\Z/p\Z$, where by $\lceil c \rceil$ for $c \in \R$ we denote the smallest integer dominating $c$. However, $\frac{2k-1}{4}+1 > \frac{2k+1}{5}$, thus $|\Z/p\Z|=p < \lceil \frac{2k-1}{4} \rceil+1$, thus we have a contradiction. \qed

 The remaining cases are when $2k+1$ is prime or a power of $3$. \\ \ \\
 \noindent
{\bf  Case 2}: $2k+1=3^l$. Let $p=3$. Suppose $\Q(e_{k-1})/\Q$ is Galois. From Proposition \ref{factor} we have 
 $$p_{k-1}(x)=m_{k-1}(x)\equiv x(x-1)^\alpha (x+1)^\beta \mod  3,$$
 where $\alpha+\beta+1=2k.$ Thus
 \begin{eqnarray*}
 P_{k-1}(q)&=&(q+q^{-1})(q+q^{-1}-1)^\alpha (q+q^{-1}+1)^\beta \cdot q^{2k} \\
 &=& (q^2+1)(q^2-q+1)^\alpha (q^2+q+1)^\beta \\
 &\equiv&  (q^2+1)(q+1)^{2\alpha} (q-1)^{2\beta} \mod 3. 
 \end{eqnarray*}
 Note that $(q^2+1)$ is irreducible $ \mod 3$. 
 Since $3| 2k+1$, 
 $P_{k-1}(1)=-2k+1=(-2k-1)+2 = 2\mod 3;$
 thus $\beta =0$. 
 On the other hand $P_{k-1}(-1)=2k+1 =0\mod 3,$ so $\alpha \neq 0$. However, we get $\alpha < 3$ by the following computation. 
 \begin{eqnarray*}
 P''_{k-1}(q) &=& 4k(4k-1)q^{4k-2}-(4k-1)(4k-2)q^{4k-3} \\
 && \quad -(4k-2)(4k-3)q^{4k-4}  -(4k-3)(4k-4)q^{4k-5} \\
 &+& \cdots \\
 & \cdots & \\
 &+& 4 \cdot 3 q^2- 3\cdot 2 q - 2\cdot 1 q^0 -1 \cdot 0, 
 \end{eqnarray*}
 thus 
 \begin{eqnarray*}
 P''_{k-1}(-1)&=& \sum_{n=1}^k \{4n(4n-1) +(4n-1)(4n-2) -(4n-2)(4n-3)  \\
 && \quad +(4n-3)(4n-4)\}  \\
 &=& \sum_{n=1}^k  (32n^2 -24 n +8) \\
 &=& 32 \cdot \frac{k(2k+1)(k+1)}{6}-24 \cdot \frac{(k+1)k}{2}+8k \\
 &=& \frac{2(2k+1)(8k-1)k}{3}+2k \\
 &\equiv& 
 \begin{cases}
 1\mod 3,  & {\text{if}} \ 2k+1=3, \\
 2k\equiv 2 \not\equiv 0 \mod 3, & {\text{if}} \ 2k+1 > 3. 
 \end{cases}
 \end{eqnarray*}
 Therefore we need $2k < 1+3$. Thus $\Q(e_{k-1})/\Q$ cannot be Galois for $k-1>1$, where $2k+1$ is a power of $3$. 
 \\ \ \\
{\bf Case 3}: $2k+1$ is a prime $\neq 3$. Let $p=2k+1$, and assume that $\Q(e_{k-1})/\Q$ is Galois.
 From Proposition \ref{factor} we have 
 $$p_{k-1}(x)=m_{k-1}(x)=x \prod_{a \in \Z/p\Z, a \neq 0} (x-a)^{\beta_a} \mod p,$$
 where $\sum_a \beta_a+1=2k.$ Thus
 \begin{eqnarray*}
 P_{k-1}(q)&=&(q+q^{-1}) \prod_a (q+q^{-1}-a)^{\beta_a} \cdot q^{2k} \\
 &=& (q^2+1) \prod_a (q^2-aq+1)^{\beta_a}    \mod p. \\
 \end{eqnarray*}
\begin{lemma} 
\label{fermat}
Let $\alpha \neq 0$ be in the algebraic closure of $\Z/p\Z=:\mathbb{F}_p$. Then 
$$\alpha + \alpha^{-1} \in \mathbb{F}_p  \Leftrightarrow \alpha^{p-1}=1 \  \mbox{\rm or} \   \alpha^{p+1}=1.$$
\end{lemma}
We postpone the proof of this lemma to the end of this subsection. If $\alpha$ is a root of $P_{k-1}(q)$, it is a root of $(q^2-bq+1)$ for some $b \in \mathbb{F}_p$; thus $\alpha+\alpha^{-1} =b \in \mathbb{F}_p$. Therefore if $\beta_a \neq 0$ and $(q^2-aq+1)$ is irreducible, $(q^2-aq+1) | q^{p-1}-1$ or $((q^2-aq+1) | q^{p+1}-1.$ Any linear factor of $P_{k-1}(q)$ divides $q^{p-1}-1$ or $q^{p+1}-1$ as well.

  On the other hand we have the following:
\begin{claim}
Let $p=2k+1 \neq 3$. Then 
\begin{itemize}
\item[(1)] $\gcd (q^{p-1}-1, P_{k-1}(q)) | (q^4-1)$. 
\item[(2)] $\gcd (q^{p+1}-1, P_{k-1}(q)) | (q^4-1)(q^3-1)$  
\end{itemize}
modulo $p$. 
\end{claim}
\noindent
{\bf Proof.} \\
(1) From the Euclidean algorithm one obtains $$\gcd (q^{p-1}-1,  Q_{k-1}(q)) | (q^4-1).$$ Since $$\gcd (q^{p-1}-1, P_{k-1}(q)) | \gcd (q^{p-1}-1, Q_{k-1}(q)), $$ we are done. 
Likewise, one obtains that $$\gcd (q^{p+1}-1, Q_{k-1}(q)) | q^6+q^5+q^4-q^2-q-1, $$ and the right hand side divides $(q^4-1)(q^3-1)$.  \qed \\ \ \\
\indent
Since $q^{p-1}-1=\prod_{0 \neq b\in \mathbb{F}_p} (q-b)$, $(q-b)$ divides $P_{k-1}(q)/(q^2+1)$ only if $b=\pm 1$. Using the same computation as in the case for $k \neq 2\mod 3$, we have $P_{k-1}(1)=(-2k-1)+2 \equiv 2 \mod p$, and $P_{k-1}(-1)=2k+1 \equiv 0\mod p$, and $P''_{k-1}(-1)=2k\equiv -1\neq 0\mod p$. (Note that $3$ is invertible in $\mathbb{F}_p$.)  Thus we have
$$P_{k-1}(q)= (q^2+1) (q+1)^2 \prod_{a \neq 0, \pm 2} (q^2-aq+1)^{\beta_a}, $$
and all the terms $(q^2-aq+1)$ appearing here are irreducible in ${\mathbb F}_p[q]$. Since they cannot divide  $q^{p-1}-1$ which is a product of linear terms, they must divide  $q^{p+1}-1$, therefore $(q^4-1)(q^3-1)$. Since $(q^4-1)(q^3-1)=(q^2+1)(q-1)(q+1)(q-1)(q^2+q+1)$, we have $\beta_a=0$ if $a\neq -1.$ 
Since Lemma \ref{multi} works for $p=2k+1>3$, we still have $\beta_a \leq 4. $  Therefore we have $\deg P_{k-1}(q)=4k \leq 12$, thus $k \leq 3$. Since $k\neq 1, 2$ by assumption, the conclusion of Proposition \ref{factor} fails for all $P_k$'s except possibly for $P_2$. For $P_2$ one may directly verify that $(q^2+q+1) {\not|} P_2(q) \mod 7$, thus Proposition \ref{factor} fails in this case as well. \qed \\

\noindent
{\bf Proof of Lemma \ref{fermat}.} \\ 
($\Rightarrow$) Suppose $\alpha + \alpha^{-1}=:m \in \mathbb{F}_p$. 
Then $\alpha$ is a root of $q^2 -mq +1 =0$. Since $m^p =m$, we have  $\alpha^{2p} - m \alpha^p +1  = (\alpha^2 -m \alpha +1)^p =0$. Thus $\alpha^p$ is also a root of  $q^2 -mq +1 =0$, and hence is equal to $\alpha$ or $\alpha^{-1}$. 
 \\
($\Leftarrow$) Suppose $\alpha^{p \pm 1}\equiv 1\mod p$. Then $\alpha^{-(p \pm 1)} \equiv 1\mod p$ as well, and $\alpha^p \equiv \alpha^{\mp 1}$.   Then $(\alpha+\alpha^{-1})^p \equiv (\alpha^p + \alpha^{-p})\equiv \alpha+\alpha^{-1}\mod p$. Therefore $\alpha+\alpha^{-1}$ is a root of $q^{p}-q=\prod_{a \in \mathbb{F}_p}  (q-a)  \equiv 0\mod p$; thus it is in $\mathbb{F}_p$.  
\subsection{The case $k \equiv 2\mod 3$}
\label{kequiv2}

We still use Proposition \ref{factor} and derive a contradiction, in essentially the same way as in the previous section. Note that  $2k+1$ cannot be be divisible by $3$ in this case. Therefore we deal with two cases: whether $2k+1$ is a prime or not. Note that $P_{k-1}(q)$ is not irreducible in this case: instead,  $P_{k-1}(q)/(q^2+q+1)$ is irreducible and it corresponds to the minimal polynomial $m_{k-1}(x)$. 
\\ \ \\
\noindent
Case 1: $2k+1$ is not a prime. \\
We take a prime $p$ so that $p | 2k+1$.  We have $p \leq \lfloor \frac{2k+1}{5} \rfloor$ as explained in the Proof in \S\ref{knequiv2}. Since $\deg m_{k-1}=2k-1$, we need at least  $ \lceil \frac{2k-2}{4} \rceil +1$ distinct elements in $\Z/p\Z$ in order for $\Q(e_{k-1})$ to be Galois by Proposition \ref{factor}. However, we still have an inequality $\frac{2k-2}{4}  +1 > \frac{2k+1}{5}$; therefore there aren't sufficiently many distinct elements in $\Z/p\Z$. 
\\ \ \\
Case 2: $2k+1$ is a prime. \\
Let $p=2k+1$. The proof is exactly the same as the previous section, except for a slight difference at the very end. We have $\deg (P_{k-1})(q)/(q^2+q+1)=4k-2 \leq 12$; thus we get the same inequality $k \leq 3$. However, by assumption $k \geq 3$ and $k=3 \not\equiv 2$.  \qed
 

\end{document}